\documentclass[12pt]{article}
\usepackage{amsmath}
\usepackage{amssymb}
\usepackage[cp1251]{inputenc}
\usepackage[russian]{babel}
\newcommand{\il}[2]{\int\limits_{#1}^{#2}}

\newcommand{\ph}{\phantom{a}}
\newcommand{\phh}{\phantom{aaa}}

\newcommand{\sist}[2]{\left\{
\begin{array}{l}
{#1}\\
\ph\\
{#2}
\end{array}
\right.}

\newcommand{\und}[2]{\hskip -20pt
\begin{array}{l}
{_{^{#1}}}\\
{^{^{#2}}}
\end{array}
}

\begin{document}
MSC 34K06
\vskip 12 pt
\centerline{ \bf  Oscillation and non-oscillation  criteria for second order}
 \centerline{ \bf  linear  non homogeneous functional-differential equations}

\vskip 20pt

\centerline{\bf G. A. Grigorian}

\vskip 20pt

\noindent
Abstract:
The Riccati equation method is used to  establish   oscillation and non-oscillation criteria for   second order linear nonhomogeneous  functional-differential equations. We show that the obtained oscillation criterion is a generalization of J. S. W. Wong's oscillation criterion for second order linear nonhomogeneous ordinary differential equations. Two examples, demonstrating the aptitude of the obtained criteria, are presented.
\vskip 20pt

\noindent
Key--Words:
Riccati equations, functional-differential equations,  oscillation, interval \linebreak oscillation,  non-oscillation.

\vskip 20pt

{\bf 1. Introduction.} Let $p(t), \ph q(t) \ph f(t), \ph q_j(t),
 \ph  r_j(t), \ph   j=\overline{1,n}$,  be real-valued  locally integrable, $\alpha_j(t),\ph  j=\overline{1,n}$ be locally measurable functions on $[t_0,+\infty)$, and let $p(t) >~ 0, \linebreak \alpha_j(t) \le t, \ph  j=\overline{1,n}, \ph t \ge t_0.$
 Consider the  second order linear functional-differential equation
$$
(p(t)\phi'(t))' + q(t) \phi'(t) + \sum\limits_{j=1}^n r_j(t)\phi(\alpha_j(t)) = f(t), \phh t \ge t_0. \eqno (1.1)
$$
Let $\theta(t)$ be a continuous function on $(-\infty, t_1],$ for some $t_1 \ge t_0$ and let $\zeta \in \mathbb{R}$.
By a Cauchy problem for Eq. (1.1)  we mean to find a continuous on $\mathbb{R}$ function $\phi(t)$,  which is  continuously differentiable on $[t_1,+\infty)$ with its absolutely continuous derivative on $[t_1,+\infty)$, satisfies (1.1) almost everywhere on $[t_1,+\infty)$ and satisfies the initial conditions $\phi(t) = \theta(t),\ph t \le t_1, \ph  \phi'(t_1) =~\zeta.$
It  should be understand  Eq. (1.1) in the sense of its equivalence to the following linear system of functional-differential equations
$$
\sist{\phi'(t) = \frac{1}{p(t)} \psi(t),}{\psi'(t) = - \sum\limits_{j=1}^n r_j(t)\phi(\alpha_j(t))-\frac{q(t)}{p(t)}\psi(t) + f(t), \phh t \ge t_0},
$$
which can be realized by involving a new dependent variable $\psi(t)$ by
$$
\psi(t) = p(t)\phi'(t), \ph t \ge t_0.
$$
This on the basis of results of work [7]
yields that for every $\zeta \in \mathbb{R}, \ph t_1\ge t_0$ and  real-valued continuous function $\theta(t)$ on $[t_1,+\infty)$ the Cauchy problem   for Eq. (1.1)   with the initial conditions   $\phi(t) = \theta(t),\ph t \le t_1, \ph  \phi'(t_1) =~\zeta$   has always a unique solution.
By a solution of Eq. (1.1) we mean any its solution of the Cauchy problem (for any $t_1 \ge t_0$).

{\bf Definition 1.1.} {\it A solution of Eq. (1.1) is called oscillatory if it has arbitrarily large zeroes.}

{\bf Definition 1.2.} {\it Eq. (1.1) is called oscillatory if its all solutions are oscillatory, otherwise it is called nonoscillatory.}

{\bf Definition 1.3.} {\it Eq. (1.1) is called oscillatory on the interval $[a,b] \subset [t_0,+\infty)$, if its every solution vanishes on $[a,b]$.}

One of the most important problem in the study of linear functional differential equations, in particular of Eq. (1.1), is that of finding conditions, providing oscillatory or nonoscillatory behavior of their solutions.
 To this problem are devoted many works (see  [1~-3, 5, 9-12] and cited works therein). One of approaches in studying of this problem is that of finding direct conditions on coefficients of the studying equation, providing its oscillation (or non-oscillation). Results in this direction have been obtained in [1, 3, 9-12]. Another approach is that of  comparing the studying equation with the another (functional-differential or ordinary differential) equation. This approach allows by means of properties of solutions of relatively simple equation to describe (to detect) wide classes of oscillatory and (or) nonoscillatory equations. Results of this type have been obtained, for example, in [1, 5, 10]. Finally,  a radical approach is that of reducing the oscillation problem for functional-differential equations to the oscillation problem for ordinary differential equations. A  result in this direction  has been obtained in [2] (see [2, Theorem 2]).

In this paper  we use the Riccati equation method  to establish  oscillation and non-oscillation criteria for   Eq.  (1,1).

 Let $d(t)> 0 , \ph r(t)$ and $g(t)$ be real-valued continuous functions on $[t_0,+\infty)$. Consider the second order linear ordinary differentia equation
$$
(d(t)\phi')' + r(t)\phi = g(t), \ph t \ge t_0. \eqno (1.2)
$$
Combining the Riccati equation method with a variational technique  J. S. W. Wong proved in [13] the following oscillation theorem

{\bf Theorem 1.1 ([8, Theorem 1])}. {\it Suppose that   for every $T\ge t_0$ there exist $T\le s_1 < t_1 \le s_2 <t_2$ such that
$$
g(t)\sist{\le 0, \ph t \in [s_1,t_1],}{\ge 0, \ph t \in [s_2,t_2].}
$$
Denote $D(s_i,t_i) \equiv \{ u \in \mathbb{C}^1[s_i,t_i] | u(t)\not \equiv 0, u(s_i) = u(t_i) =0\}, \ph i=1,2.$ If there exists $u\in  D(s_i,t_i)$ such that
$$
\il{s_i}{t_i}(r(\tau)u(\tau)^2 - d(\tau) u'(\tau)^2) d \tau \ge 0
$$
for $i=1,2$, then Eq. (1.2) is oscillatory.}

We show that the obtained below oscillation criterion is a generalization of
 Theorem~ 1.1.

\vskip 20pt

{\bf 2. Auxiliary propositions}. Let   $a(t), b(t), c(t), a_1(t), b_1(t), c_1(t)$ be real-valued locally integrable functions on  $[t_0, +\infty)$. Consider the Riccati equations
$$
y'(t) + a(t) y^2(t) + b(t) y(t) + c(t) = 0, \eqno (2.1)
$$
$$
y'(t) + a_1(t) y^2(t) + b_1(t) y(t) +  c_1(t) = 0, \eqno (2.2)
$$
$t\ge t_0$
and the differential inequalities
$$
\eta'(t) + a(t) \eta^2(t) + b(t) \eta(t) + c(t) \ge 0, \eqno (2.3)
$$
$$
\eta'(t) + a_1(t) \eta^2(t) + b_1(t) \eta(t) + c_1(t) \ge 0,  \eqno (2.4)
$$
$t\ge t_0.$
By a solution of Eq. (2.1) (Eq. (2.2), inequalities (2.3), (2.4)) on $[t_1,t_2)\ph (t_0 \le t_1 < t_2\le +\infty)$ we mean an absolutely continuous function on $[t_1,t_2)$, satisfying (2.1) ( (2.2), (2.3), (2.4)) almost everywhere on $[t_1,t_2)$. Since the function $a(t)y^2 + b(t) y + c(t) \ph  (a_1(t)y^2 + b_1(t) y + c_1(t))$ satisfies the Caratheodory condition on $[t_0,+\infty)\times \mathbb{R}$, for every $t_1\ge t_0, \ph \gamma \in \mathbb{R}$ there exists $t_2 > t_1$ such that Eq. (2.1) (Eq. (2.2)) has a solution $y(t)$  on $[t_1,t_2)$ with $y(t_1) =\gamma$.
Note, that every solution of Eq. (2.1) ( (2.2) ) is a solution of the inequality  (2.3)\phantom{a} (  (2.4) ). Note also, that for  $a(t) \ge 0 \phantom{a}(a_1(t) \ge 0), \phantom{a}t \ge t_0$, the real-valued  solutions of the equation  $\eta'(t) + b(t) \eta(t) + c(t) = 0$ \phantom{a} ($\eta'(t) + b_1(t) \eta(t) + c_1(t) = 0$)  are solutions of the inequality (2.3) \ph ( (2.4) ).
Therefore for  $a(t) \ge 0 \phantom{a} (a_1(t) \ge 0),\phantom{a} t \ge t_0$,  the inequality   (2.3) (  (2.4) ) has a solution on $[t_0,+\infty)$, satisfying any initial real-valued condition. In the sequel we will assume, that the solutions of considering equations and inequalities are real-valued.

{\bf Theorem 2.1}. {\it Let $y_0(t)$  be a solution of Eq. (2.1) on  $[t_1, t_2)$, and  $\eta_0(t)$, $\eta_1(t)$  be solutions of inequalities (2.3)   and  (2.4) with  $\eta_0(t_1) \ge y_0(t_1), \phantom{a}\eta_1(t_1) \ge y_0(t_1)$ respectively, and let $a_1(t) \ge 0, \phantom{a} \lambda - y_0(t_1) +$\\
$
+ \int\limits_{t_1}^t\exp\biggl\{\int\limits_{t_1}^\tau[a_1(\xi)(\eta_0(\xi) +\eta_1(\xi)) + b_1(\xi)] d \xi\biggr\} \times
$
$$
\phantom{aaaaaaaaaa}\times [(a(\tau) - a_1(\tau)) y_0^2(\tau) + (b(\tau) - b_1(\tau)) y_0(\tau)  + c(\tau) - c_1(\tau)] d\tau \ge 0, \phantom{a}t\in [t_1,t_2),
$$
for some $\lambda \in [y_0(t_1),\eta_1(t_1)]$.
Then Eq. (2.2) has a solution  $y_1(t)$ on  $[t_1,t_2)$  with  $y_1(t_1) \ge\linebreak \ge y_0(t_1)$, moreover  $y_1(t) \ge   y_0(t),\phantom{a} t\in [t_1,t_2)$.}

Proof. By analogy of the proof of Theorem  3.1 from [4].

Consider the equation
$$
y'(t) + \frac{1}{p(t)} y^2(t) + \frac{q(t)}{p(t)} y(t)
 -\frac{f(t)}{\lambda}\exp\bigg\{ - \int\limits_{t_1}^t\frac{y(\tau)}{p(\tau)}d\tau\biggr\}  + \phantom{aaaaaaaaaaaaaaaaaaaaaaa}
$$
$$
 \phantom{aaaaaaaa}+ \sum\limits_{j=1}^n r_j(\tau)\exp\biggl\{-\il{\alpha_j(t)}{t}\frac{y(\tau)}{p(\tau)}d\tau\biggr\} = 0, \ph  t\ge t_1 \phantom{a} (\ge t_0), \phantom{a} \lambda = const \ne 0. \eqno (2.5)
 $$

By a solution of  this equation on $[t_1,t_2) \ph (t_0 \le t_1 < t_2 \le +\infty)$ we mean a continuous function $y(t)$ on $(-\infty,t_2)$, which is absolutely continuously  on $[t_1,t_2)$ and satisfies (2.5) almost everywhere on $[t_1,t_2)$.

Let  $t_0 \le t_1 < t_2 \le +\infty$. Denote $T(t_1,t_2)
\equiv \min\{t_1,\min\limits_{1\le j \le n}\{\inf\limits_{t\in [t_1,t_2)}\alpha_j(t)
\}.
$

Let  $\phi_0(t)$ be a solution of Eq.  (1.1) on  $[t_1,t_2)$,
and let  $\phi_0(t) \ne 0,\phantom{a} t\in  [T(t_1,t_2),t_2)$.  It is easy to show, that
$$
y_0(t)\equiv\sist{\frac{p(t)\phi'_0(t)}{\phi_0(t)}, \ph t\in [t_1,t_2),}{ \ph \frac{p(t_1)\phi'_0(t_1)}{\phi_0(t_1)}, \ph  t \le t_1} \eqno (2.6)
$$
is a solution of Eq. (2.5) on  $[t_1,t_2)$, where $\lambda = \phi_0(t_1)$.

{\bf Definition 2.1.} {\it An interval $[t_1,t_2), \ph (t_0 \le t_1 <t_2 \le +\infty)$ is called the maximum existence interval for a solution $y(t)$ of Eq. (2.5) if it exists on $[t_1,t_2)$ and cannot be continued to the right from $t_2$ as a solution of Eq. (2.5).}

{\bf Lemma 2.1.} {\it Let $y(t)$ be a solution of Eq. (2.5) on $[t_1,t_2)$ with $t_2 < +\infty$. If the function $F(t) \equiv \il{t_1}{t}\frac{y(\tau)}{p(\tau)} d \tau, \ph t_1\le t < t_2$ is bounded from below, then $[t_1,t_2)$ cannot be the maximum existence interval for $y(t)$.}

Proof. By (2.6) $\phi(t) \equiv \exp\bigl\{F(t)\}, \ph t_1\le t < t_2$ is a solution to Eq. (1.1). From here and from the condition of the lemma it follows that $\phi(t) \ne 0, \ph t \in [t_1, t_3)$ for some $t_3 > t_2$. Then by (2.6) the function $\frac{p(t)\phi'(t)}{p(t)}, \ph t_1\le t < t_3$ is a solution of Eq (2.5) on $[t_1,t_3)$, which is a continuation of $y(t)$. Therefore $[y_1,t_2)$ is not the maximum existence interval for $y(t)$. The lemma is proved.

Let $r_{1,j}(t), \ph j= \overline{1,n}, \ph t \ge t_0$ be real-valued continuous functions on $[t_0,+\infty)$.
Consider the equation
$$
y'(t) + \frac{1}{p(t)} y^2(t) - \frac{q(t)}{p(t)}y(t)
+ \sum\limits_{j=1}^n r_{1,j}(\tau)\exp\biggl\{-\il{\alpha_j(t)}{t}\frac{y(\tau)}{p(\tau)}d\tau\biggr\} = 0, \ph t\ge t_1 \phantom{a} (\ge t_0). \eqno (2.7)
 $$

{\bf Lemma 2.2}. {\it Let $y_1(t)$ be a solution of  Eq. (2.7) on   $[t_1,t_2)$,
and let the following conditions be satisfied.

\noindent
1) $r_{1,j}(t) \ge r_j(t),\ph t \in [t_1,t_2),\ph j=\overline{1,n}$.

\noindent
2) If $r_j(t) < 0,$ then $r_{1,j}(t) \ge 0$ or $\alpha_j(t) = t$ for every $t \in [t_1,t_2)\ph j=\overline{1,n}$.

\noindent
3) $ \frac{f(t)}{\lambda}\ge 0, \ph t \in [t_1,t_2)$.

\noindent
Then for every real-valued continuous function $\gamma(t) \ge y_1(t), \ph t\le t_1, \ph \gamma(t_1) > y_1(t_1)$ Eq.~ (2.5) has a solution $y(t)$ on $[t_1,t_2)$ with $y(t) =~ \gamma(t), \ph t \le t_1$ and
$$
y(t) > y_1(t), \phh t \in [t_1,t_2). \eqno (2.8)
$$
}

Proof. Let $y(t)$  be a  solution of Eq. (2.5)  with $y(t) = \gamma(t), \ph t \le t_1$, and let $[t_1,t_3)$ be its maximum existence interval. Suppose $t_3 < t_2$. Show that
$$
y(t) > y_1(t), \phh t \in [t_1,t_3). \eqno (2.9)
$$
Set:
$$
R(t) \equiv \sum\limits_{j=1}^n r_j(\tau)\exp\biggl\{-\il{\alpha_j(t)}{t}\frac{y(\tau)}{p(\tau)}d\tau\biggr\}  - \frac{f(t)}{\lambda}\exp\bigg\{ - \int\limits_{t_1}^t\frac{y(\tau)}{p(\tau)}d\tau\biggr\}, \phantom{aaaaaaaaaaaaaaaaaaaaaaaa}
$$
$$
\phantom{aaaaaaaaaaaaaaaaaaaaaaaaaa}R_1(t) \equiv \sum\limits_{j=1}^n r_{1,j}(\tau)\exp\biggl\{-\il{\beta_j(t)}{t}\frac{y_1(\tau)}{p(\tau)}d\tau\biggr\}, \ph t \in [t_1,t_3).
$$
Suppose (2.9) is not true. Since $y(t_1) > y_1(t_1)$ there exists $t_4 \in (t_1,t_3)$ such that
$$
y(t) > y_1(t), \ph t \in [t_1,t_4), \ph \mbox{and} \ph y(t_4) = y_1(t_4). \eqno (2.10)
$$
Show that
$$
R(t) \le R_1(t), \phh t \in [t_1,t_4). \eqno (2.11)
$$
Set
$$
\Delta_j(t) \equiv r_{1,j}(t)\exp\biggl\{-\il{\alpha_j(t)}{t}\frac{y_1(\tau)}{p(\tau)}d\tau\biggr\} - r_j(t)\exp\biggl\{-\il{\alpha_j(t)}{t}\frac{y(\tau)}{p(\tau)}d\tau\biggr\}, \ph t \in [t_1,t_3).
$$
Obviously if $r_{1,j}(t) \ge 0, \ph r_j(t) < 0$, then $\Delta_j(t) \ge 0.$ If $r_j(t) < 0$ and $\alpha_j(t) = t$, them by the condition 1) we have $\Delta = r_{1,j}(t) - r_j(t) \ge 0.$ Assume $r_j(t) \ge 0$. Then
$$
\Delta_j(t) = [r_{1,j}(t) - r_j(t)] \exp\biggl\{-\il{\alpha_j(t)}{t}\frac{y_1(\tau)}{p(\tau)}d\tau\biggr\} + r_j(t)\exp\biggl\{-\il{\alpha_j(t)}{t}\frac{y(\tau)}{p(\tau)}d\tau\biggr\}\times \phantom{aaaaaaaaaaaaaaaaaaaa}
$$
$$
\phantom{aaaaaaaaaaaaaaaaaaaaaaaaaaaaaaaaaaaaaaaa}\times\biggl[1 -\exp\biggl\{-\il{\alpha_j(t)}{t}\frac{[y_1(\tau)- y(\tau)]}{p(\tau)}d\tau\biggr\}\biggr]
$$
This together with the condition 1) and (2.10) implies that $\Delta_j(t) \ge 0$. Hence, (2.11) is valid.

Consider the Riccati equations
$$
y' + \frac{1}{p(t)} y^2 + q(t) y + R(t) = 0, \ph t \in [t_1,t_3), \eqno (2.12)
$$
$$
y' + \frac{1}{p(t)} y^2 + q(t) y + R_1(t) = 0, \ph t \in [t_1,t_3). \eqno (2.13)
$$
Note that $y(t)$ and $y_1(t)$ are solutions of the equations (2.12) and (2.13) respectively on $[t_1,t_3)$. Then since $y(t_1) > y_1(t_1)$  by Theorem 2.1 from (2.11) it follows that  $y(t_4) > y_1(t_4)$, which contradicts (2.10). The obtained contradiction proves (2.9). It follows from (2.9) that the function $F(t) \equiv \il{t_1}{t}\frac{y(t)}{p(t)}, \ph t \in [t_1,t_3)$ is bounded from below on $[t_1,t_3)$. In virtue of Lemma 2.1 from here it follows that $[t_1, t_3)$ is not the maximum existence interval for $y(t)$, which contradicts our supposition. The obtained contradiction shows that $y(t)$ exists on $[t_1,t_2)$ and the inequality (2.8) is valid. If $t_3 \ge t_2$ then the existence $y(t)$ on $[t_1,t_2)$ is clear and the inequality (2.8) follows from the already proven that for any $t_5 \in (t_1,t_2)$ the inequality $y(t) > y_1(t), \ph t \in [t_1,t_5)$ is valid. The lemma is proved.

By analogy of the proof of Lemma 2.2 can be proved.

{\bf Lemma 2.3}. {\it Let $y(t)$ be a solution of  Eq. (2.5) on   $[t_1,t_2]$,
and let the following conditions be satisfied.

\noindent
1') $r_{1,j}(t) \le r_j(t),\ph t \in [t_1,t_2),\ph j=\overline{1,n}$.

\noindent
2') If $r_{1,j}(t) < 0,$ then $r_{j}(t) \ge 0$ or $\alpha_j(t) = t$ for every $t \in [t_1,t_2]\ph j=\overline{1,n}$.

\noindent
3') $ \frac{f(t)}{\lambda}\le 0, \ph t \in [t_1,t_2]$.

\noindent
Then for every real-valued continuous function $\gamma(t) \ge y(t), \ph t\le t_1, \ph \gamma(t_1) > y(t_1)$ Eq.~ (2.7) has a solution $y_1(t)$ on $[t_1,t_2]$ with $y_1(t) =~ \gamma(t), \ph t \le t_1$ and
$$
y_1(t) > y(t), \phh t \in [t_1,t_2].
$$
}
\phantom{aaaaaaaaaaaaaaaaaaaaaaaaaaaaaaaaaaaaaaaaaaaaaaaaaaaaaaaaaaaaaaaaaaaaaaaaaaa} $\blacksquare$

Let $\beta_j(t), \ph j=\overline{1,n}$ be locally measurable functions on $[t_0,+\infty)$ and let  $t_0 \le t_1 < t_2 \le t_3 < t_4 < + \infty$.
Set: $\omega_+\equiv \{k\in\Omega | \beta_k(t) \ge t, \ph t\in [t_1,t_2]\};$
$\ph\omega_1^-\equiv\{k\in \Omega : t_1\le \beta_k(t) \le t, \ph t\in [t_1,t_2]\}; \ph \omega_2^-\equiv \{k\in \Omega : \beta_k(t) \le t, \ph t\in[t_3, t_4]\}; $

\noindent
$
 T_1\equiv\inf_{\und{t\in [t_1,t_4]}{k\in \Omega}} \beta_k(t);\phh
 \ph T_2\equiv \sup_{ \und{t\in [t_1,t_4]}{k\in \Omega}} \beta_k(t); \phh t_2^+ \equiv \sup_{ \und{t\in [t_1,t_2]}{k\in \omega_+}} \beta_k(t); \phh t_3^- \equiv \inf_{\und{t\in [t_3,t_4]}{k\in \omega_2^-}} \beta_k(t).
$

Consider the equations

$$
(p(t)\phi'(t))' +  \sum\limits_{k=1}^n r_k(t)\phi(\beta_k(t)) = f(t), \phh t \ge t_0. \eqno (2.14)
$$
$$
(p(t)\phi'(t))' + \Biggl[\sum\limits_{k\in\omega_+}r_k(t)\exp\Bigl\{\il{t}{\beta_k(t)} \frac{d \tau}{p(\tau)}\il{\tau}{t_2}\Bigl(\sum\limits_{k\in\omega_+}r_k(s)\Bigr) d s\Bigr\} + \phantom{aaaaaaaaaaaaaaaaaaaaaaaaaaaa}
$$
$$
\phantom{aaaaaaaaaaaaaaaaaaa}+\sum\limits_{k\in \omega_1^-} r_k(t) \frac{\il{t_1}{\beta_k(t)}\frac{d \tau}{p(\tau)} + \varepsilon}{\il{t_1}{t}\frac{d \tau}{p(\tau)}+ \varepsilon}\Biggr] \phi(t) = 0, \phh t \in [t_1, t_2], \phh \varepsilon > 0; \eqno (2.15)
$$
$$
(p(t)\phi'(t))' + \Bigl[\sum\limits_{k\in\omega_2^-} r_k(t)\Bigr] \phi(t) = 0, \phh t\in [t_3,t_4]. \eqno (2.16)
$$

{\bf Theorem 2.2. ([5, Theorem 2.6])} {\it Let $\omega_+ \cup \omega_1^- \ne \emptyset, \phantom{a} \omega_2^- \ne \emptyset$ and let the following conditions be satisfied:

\noindent
a) \ph $r_k(t) \ge 0, \ph t\in [T_1, T_2], \ph k \in\Omega;$

\noindent
b) \ph $t_2^+ \le t_3^-$;

\noindent
c) \ph for some $\varepsilon_0 > 0$  Eq. (2.15) is oscillatory on $[t_1,t_2]$  for all  \ph $\varepsilon \in (0,\varepsilon_0)$;

\noindent
d) \ph Eq. (2.16) is oscillatory on  $[t_3, t_4]$.

\noindent
Then Eq. (2.14) is oscillatory on $[T_1,T_2]$.}

{\bf 3. Oscillation and non-oscillation criteria}.

Along with Eq. (1.1) consider the equation
$$
(p(t)\phi'(t))' + q(t)\phi'(t) + \sum\limits_{j=1}^n r_{1,j}(t)\phi(\alpha_j(t)) =0, \phh t\ge t_0. \eqno (3.1)
$$

{\bf Theorem 3.1.} {\it Let the following  conditions be satisfied.

\noindent
1) $r_{1,j}(t) \ge r_j(t),\ph t \ge t_0, \ph j=\overline{1,n}$.

\noindent
2) If $r_j(t) < 0,$ then $r_{1,j}(t) \ge 0$ or $\alpha_j(t) = t$ for every $t \in [t_1,t_2)\ph j=\overline{1,n}.$

\noindent
4) $f(t) \ge 0, \ph t \ge t_0$.

\noindent
5) $\lim\limits_{t \to +\infty}\alpha_j(t) = +\infty, \ph j=\overline{1,n}$.

 If   Eq. (3.1) is nonoscillatory, then Eq. (1.1) is also nonoscillatory.}

Proof. Let Eq. (3.1) be nonoscillatory. Then it has a solution $\phi_1(t)$ such that $\phi_1(t) \ne~ 0, \linebreak t \ge t_1,$ for some $t_1\ge t_0$, and, therefore,  the condition 5) implies that
$$
\phi(\alpha_j(t)) \ne 0, \phh t \ge T, \phh j=\overline{1,n}.
$$
for some $T > t_1$.
By virtue of (2.6) from here it follows that Eq. (2.7) has a solution $y_1(t)\equiv \frac{p(t)\phi'(t)}{p(t)}$ on $[T,+\infty)$.

 Then by Lemma 2.2 from the conditions 1), 2) and 4)   it follows that Eq. (2.5) has a solution $y(t)$ on $[T,+\infty)$. Hence by (2.6)
$$
\phi(t) \equiv \exp\biggl\{\il{T}{t}\frac{y(\tau)}{p(\tau)} d \tau\biggr\}, \ph t \in \mathbb{R}
$$
is a nonoscillatory solution of Eq. (1.1) on $[T,+\infty).$
Therefore Eq. (1.1) is nonoscillatory. The theorem is proved.

Set $r_j^+(t) \equiv \max\{0,r_j(t)\}, \ph t\ge t_0, \ph j=\overline{1,n}$. Consider the equation
$$
(p(t)\phi'(t))' + q(t)\phi'(t) + \sum\limits_{j=1}^n r_j^+(t)\phi(\alpha_j(t)) =0, \phh t\ge t_0. \eqno (3.2)
$$
Since $r^+_j(t)\ge 0, \ph t \ge t_0, \ph j=\overline{1,n},$ from Theorem 3.1 we obtain immediately

{\bf Corollary 3.1}. {\it Let the  conditions 4) and 5) of Theorem 3.1 be satisfied.
If Eq. (3.2) is nonoscillatory, then Eq. (1.1) is also nonoscillatory.
}

\phantom{aaaaaaaaaaaaaaaaaaaaaaaaaaaaaaaaaaaaaaaaaaaaaaaaaaaaaaaaaaaaaaaaaaaaaaa} $\blacksquare$

{\bf Example 3.1.} {\it Consider the equation
$$
\phi''(t) + \sin^2 t \hskip 2.5pt\phi(\alpha_1(t)) + \cos^2 t\hskip 2.5pt \phi(\alpha_2(t)) - \phi(t)  = cos(\sin (\ln(1+t))), \ph t \ge 0, \eqno (3.3)
$$
$\lim\limits_{t \to + \infty}\alpha_j(t) = + \infty, \ph j=1,2.$
For establishing  non-oscillation of this equation we compare it with the following homogeneous one
$$
\phi''(t) + \sin^2 t \hskip 2.5pt\phi(\alpha_1(t)) + \cos^2 t \hskip 2.5pt\phi(\alpha_2(t)) - \phi(t) = 0, \ph t \ge 0, \eqno (3.4)
$$
It is not difficult to verify that for $r_{1,1}(t) = r_1(t) = \sin^2 t, \ph r_{1,2}(t) = r_2(t) = \cos^2 t, \ph r_{1,3}(t) = r_3(t) = -1, \ph f(t) = cos(\sin (\ln(1+t))), \ph t \ge 0$ the conditions 1) - 5) of Theorem 3.1 for Eqs. (3.3) and (3.4) are satisfied. Obviously, $\phi(t) \equiv 1$ is a nonoscillating solution for Eq. (3.4). Then by Theorem 3.1 Eq. (3.3) is nonoscillatory.}

{\bf Theorem 3.2.} {\it Let the following conditions be satisfied.

\noindent
I) $r_j(t) \ge r_{1,j}(t), \ph t \ge t_0 \ph j=\overline{1,n}.$

\noindent
II) If $r_{1,j}(t) < 0,$ then $r_j(t) \ge 0,$ or $\alpha_j(t) = t$, for every $t\ge t_0 \ph j=\overline{1,n}$.

\noindent
III) $\lim\limits_{t \to +\infty} \alpha_j(t) = +\infty, \ph j=\overline{1,n}$,

\noindent
IV) for every $T \ge t_0$  there exists  $T < s_1 < t_1 \le s_2 < t_2$ such that
 $$f(t) \sist{\le 0, \ph t \in [s_1,t_1],}{\ge 0, \ph t\in [s_2,t_2].}$$

\noindent
V) Eq. (3.1)
is oscillatory on the intervals $[s_k,t_k], \ph k=1,2$.

\noindent
Then Eq. (1.1) is oscillatory.}

Proof. Suppose Eq. (1.1) is not oscillatory. Then there exists its a solution $\phi(t)$ and $t_1\ge t_0$ such that $\phi(t) \ne 0, \ph t \ge t_1$. Therefore it follows from the condition III) that
$$
\phi(\alpha_j(t)) \ne 0, \ph t \ge T, \ph j=\overline{1,n} \eqno (3.5)
$$
for some $T > t_1.$ Consider the equations
$$
y' + \frac{1}{p(t)} y^2 + \frac{q(t)}{p(t)} y + \sum\limits_{j=1}^nr_j(t)\exp\biggl\{-\il{\alpha_j(t)}{t}\frac{y(\tau)}{p(\tau)} d\tau\biggr\} = \frac{f(t)}{\phi(T)}, \phh t \ge T, \eqno (3.6)
$$
$$
y' + \frac{1}{p(t)} y^2 + \frac{q(t)}{p(t)} y + \sum\limits_{j=1}^nr_{1,j}(t)\exp\biggl\{-\il{\alpha_j(t)}{t}\frac{y(\tau)}{p(\tau)} d\tau\biggr\} = 0, \phh t \ge T, \eqno (3.7)
$$
By Lemma 2.1 from (3.5) it follows, that Eq. (3.6) has a solution  on $[T,+\infty)$. Suppose $\phi(t) > 0 (< 0) \ph t \ge t_1$. Then $\frac{f(t)}{\phi(T)} \le 0, \ph t \in [s_1,t_1] \ph (t \in [s_2,t_2])$ By Lemma 2.3 from here and from the conditions I) and III) it follows that Eq. (3.7) has a solution $y_1(t)$ on $[s_1,t_1]$ (on $[s_2,t_2]$). By (2.6) from here it follows that  $\phi_1(t)\equiv \exp\biggl\{\il{T}{t}\frac{y_1(\tau)}{p(\tau)}\biggr\}, \ph t \le t_1 \ph (t \le t_2)$ is a nonoscillating on $[s_1,t_1]$ (on $[s_2,t_2]$) solution of Eq. (3.1).  Hence,
Eq. (3.1) is not oscillatory on $[s_1,t_1]$ (on $[s_2,t_2]),$ which contradicts the condition V). The obtained contradiction completes the proof of the theorem.

Consider the equations
$$
(p(t)\phi'(t))' + q(t)\phi'(t) + \sum\limits_{j=1}^n r_j^+(t)\phi(\alpha_j(t)) =0, \phh t\ge t_0. \eqno (3.8)
$$
$$
(p(t)\phi'(t))' + q(t)\phi'(t) + \sum\limits_{j=1}^n r_j^+(t)\phi(\alpha_j(t)) =f(t), \phh t\ge t_0. \eqno (3.9)
$$
From Theorem 3.2 we obtain immediately

{\bf Corollary 3.2.} {\it Let the conditions III) and IV) of Theorem 3.2 be satisfied.

\noindent
If Eq. (3.8) is oscillatory on the intervals $[s_i,t_i], \ph i=1,2$, then Eq. (3.9) is oscillatory.}

\phantom{aaaaaaaaaaaaaaaaaaaaaaaaaaaaaaaaaaaaaaaaaaaaaaaaaaaaaaaaaaaaaaaaaaaaaaa} $\blacksquare$

In the case $n=1, \ph \alpha_1(t) = t$ and continuous $r_{1,1}(t) = r_1(t) \ph t \ge t_0$ Theorem 3.2 is the same as Corollary 3.1 from [6]. It was shown in [6] that Corollary 3.1 of [6] is a generalization of Theorem 1.1. Therefore Theorem 3.2 is another generalization of Theorem 1.1.

{\bf Example 3.2.} {\it Let $c_k(t), \ph h_k(t) \ph k=\overline{1,n}$ be real-valued locally integrable functions on $[0,+\infty).$ Consider the equation
$$
\phi''(t) + \sum\limits_{k=1}^n c_k(t) \phi(t - h_k(t)) = \sin (t/3), \ph t \ge 0. \eqno (3.10)
$$
Assume $c_k(t) \equiv 0, \ph t \in [3\pi l, 3\pi l +1], \phh k=\overline{1,n} \phh c_k(t) \ge 0, \ph \sum\limits_{k=1}^n c_k(t) \ge 2, \linebreak t \in [3\pi l +1,3\pi (l+1)], \ph l=0,1,2, \dots, \ph 0\le h_k(t) \le \frac{1}{2}, \ph t \ge 0, \ph k=\overline{1,n}$.
Set: $t_{1,l}\equiv 3\pi l + \frac{1}{2}, \ph t_{2,l} \equiv (3 l + 2)\pi + \frac{1}{2}, \ph t_{3,l} \equiv (3 l + 2)\pi + 1, \ph t_{4.l}\equiv 3\pi (l+1),\ph l=0,1,\dots$. Then
$$
 T_{1,l}\equiv\inf_{\und{t\in [t_{1,l},t_{4,l}]}{k\in \Omega}}(t - h_k(t)) \ge 3\pi l;\phh
 \ph T_{2,l}\equiv \sup_{ \und{t\in [t_{1,l},t_{4,l}]}{k\in \Omega}} (t -h_k(t)) \le 3\pi (l+1), \ph l=0,1,\dots.  \eqno (3.11)
$$

It is not difficult to verify that  for $t_s = t_{s,l}, \ph s=\overline{1,4}$    the conditions of Theorem 2.2 for the  homogeneous equation
$$
\phi''(t) + \sum\limits_{k=1}^n c_k(t) \phi(t - h_k(t)) = 0, \ph t \ge 0 \eqno (3.12)
$$
are satisfied for every $l=0,1,\dots$. Then by Theorem 2.2  from (3.11) it follows that the last equation is oscillatory on each of intervals $[3\pi l,3\pi(l+1)], \ph l=0,1,2, \dots$. It follows from here that the conditions of Theorem 3.2 for Eqs. (3.10) and (3.12) are satisfied.
Therefore, Eq. (3.10) is oscillatory.}

\vskip 20pt

\centerline{\bf References}

\vskip 20pt

\noindent
1. L. Berezansky, E. Braverman. Some oscillation problems for a second order linear delay \linebreak \phantom{a} differential equations. Mathematical analysis and applications, 220, pp. 719 - 740,  1998,

\noindent
2. L. Berezansky, E. Braverman. Oscillation of a Second Order Delay Differential Equations \linebreak \phantom{a}  with Middle Therm. Applied Mathematics Letters, 13 (2000) 21-25.

\noindent
3 J. Dzurina, Oscillation Theorems for Decond order Advanced Neutral Differential \linebreak \phantom{aa}  Equations. Tatra Mt. Math. Publ., 48 (2011), pp. 61-71.

\noindent
4. G. A. Grigorian, On two comparison tests for second-order linear  ordinary differential\linebreak \phantom{a} equations (Russian) Differ. Uravn. 47 (2011), no. 9, 1225 - 1240; translation in Differ.\linebreak \phantom{a} Equ. 47 (2011), no. 9 1237 - 1252, 34C10.

\noindent
5. G. A. Grigorian, Oscillation Criteria for the Second Order Linear Functional-Differential \linebreak \phantom{a} Equations With Locally Integrable Coefficients. Sarajevo J. Math, Vol.14 (27),\linebreak \phantom{a}  No.1, (2018), pp. 71–86.

\noindent
6. G. A. Grigorian,    Oscillation and non-oscillation criteria for linear nonhomogeneous \linebreak \phantom{a} systems of two first-order ordinary differential equations  J. Math.Anal.Appl. 507 \linebreak \phantom{a} (2022) 125734,

\noindent
7. G. A. Grigorian,  Cauchy problem for quasilinear systems of functional differential \linebreak \phantom{aa} equations.  Sarajevo J. Math. Accepted for publication.

\noindent
8. Kong, Qingkai,  Pa$\check s$i$\acute c$, Mervan. Second-order differential equations: some significant \linebreak \phantom{aa} results due to James S. W. Wong. Differ. Equ. Appl. 6 (2014), no. 1, 99--163.

\noindent
9. T. Li, Yu V. Rogovchenko and Ch. Zheng, Oscillation of Second-Order Neutral Differential \linebreak \phantom{aa}  Equations. Funkcialaj Ekvacioj, 56 (2013), pp. 111-120.

\noindent
10. J. Ohrishka. Oscillation of second-order linear delay differential equations. Cent. Eur.\linebreak \phantom{aa} J. Math. N$^\circ$ 6 (3),  2008, pp. 439 - 452.

\noindent
11.  Z. Oplustil, J. Sremir. Some oscillation criteria for the second-order linear delay\linebreak \phantom{aa} differential equation. Mathematica Bohemica, vol. 136,  N$^\circ$ 2, pp. 195 - 204,  20011.

\noindent
12. K. Takasi, V. Maric. On a Class of Functional Differential Equations, Having Slowly \linebreak \phantom{aa}  Varying Solutions. Publications De L'Institut Mathematique. Novelle Serie, tome 80\linebreak \phantom{aa}  (94) (2006) pp. 207-217.

\noindent
13. Wong, James S. W. Oscillation criteria for a forced second-order linear differential \linebreak \phantom{aa}  equation. J. Math. Anal. Appl. 231 (1999), no. 1, 235--240.

\end{document}